 \documentclass[11pt]{article}

 \usepackage{mathrsfs}
 \usepackage{amssymb}
 \usepackage[dvips]{color}

 \newtheorem{lemma}{Lemma}[section]
 \newtheorem{theorem}{Theorem}[section]
 \newtheorem{corollary}{Corollary}[section]

 \def\blemma{\begin{lemma}\sl{}\def\elemma{\end{lemma}}}
 \def\btheorem{\begin{theorem}\sl{}\def\etheorem{\end{theorem}}}

 \def\beqlb{\begin{eqnarray}}\def\eeqlb{\end{eqnarray}}
 \def\beqnn{\begin{eqnarray*}}\def\eeqnn{\end{eqnarray*}}

 \def\mcr{\mathscr}
 \def\mbb{\mathbb}

 \def\qed{\hfill\mbox{$\Box$}\medskip}

 \def\<{\langle}\def\>{\rangle}

 \def\supp{{\mbox{\rm supp}}}

 \def\bfC{\mbox{\boldmath $C$}}\def\bfE{\mbox{\boldmath $E$}}
 \def\bfP{\mbox{\boldmath $P$}}\def\bfQ{\mbox{\boldmath $Q$}}

 \def\L{{\mcr L}}

 \def\itDelta{{\it\Delta}}\def\itGamma{{\it\Gamma}}
 \def\itOmega{{\it\Omega}}
 \def\itPhi{{\it\Phi}}

 \def\dfR{\mbb{R}}

 \def\lin{----------------}

\begin{document}

\noindent {Published in: {\it Journal of Theoretical Probability}
{\bf 17} (2004), 3: 673--692.}

\bigskip\bigskip

\noindent {\LARGE\bf Superprocesses with Coalescing Brownian}

\medskip

\noindent {\LARGE\bf Spatial Motion as Large-Scale Limits}

\bigskip\bigskip

\noindent Donald A. Dawson,\footnote{School of Mathematics and
Statistics, Carleton University, 1125 Colonel By Drive, Ottawa,
Canada K1S 5B6. E-mail: {\tt ddawson@math.carleton.ca}}
 \ \
Zenghu Li,\footnote{Department of Mathematics, Beijing Normal
University, Beijing 100875, People's Republic of China. E-mail:
{\tt lizh@email.bnu.edu.cn}}
 \ \ and \ \
Xiaowen Zhou\footnote{Department of Mathematics and Statistics,
Concordia University, 7141 Sherbrooke Street West, Montreal,
Canada, H4B 1R6. E-mail: {\tt zhou@alcor.concordia.ca}}

\bigskip\bigskip

{\narrower{\narrower

\centerline{\bf\lin\lin\lin\lin\lin\lin}

\medskip

\noindent {\bf Abstract.} A superprocess with coalescing spatial
motion is constructed in terms of one-dimensional excursions.
Based on this construction, it is proved that the superprocess is
purely atomic and arises as scaling limit of a special form of the
superprocess with dependent spatial motion studied in Dawson {\it
et al.} (2001) and Wang (1997, 1998).

\medskip

\centerline{\bf\lin\lin\lin\lin\lin\lin}

\bigskip

\noindent {\it Mathematics Subject Classifications (2000):}
Primary 60J80; Secondary 60G57

\bigskip

\noindent {\it Key words:} superprocess, coalescing spatial
motion, excursions, scaling limit, Poisson random measure.

\par}\par}

\bigskip\bigskip


\section{Introduction}

\setcounter{equation}{0}

Large scale limits of interacting particle systems and
measure-valued processes have been studied by many authors; see,
e.g., Bojdecki and Gorostiza (1986), Cox {\it et al.} (2000),
Dawson (1977), Dawson and Fleischmann (1988), Durrett and Perkins
(1999), Hara and Slade (2000a,b). In particular, Dawson and
Fleischmann (1988) investigated the scaling limit of a class of
critical space-time branching models, giving a precise description
of the growth of large clumps at spatially rare sites in low
dimensions. They showed that a space-time-mass scaling limit
exists and is a measure-valued branching process without
migration. The clumps are located at Poissonian points and their
sizes evolve according to continuous-state branching processes.
Durrett and Perkins (1999) proved that suitably rescaled contact
processes converge to super-Brownian motion in two or more
dimensions. Cox {\it et al.} (2000) proved convergence of some
rescaled voter models to super-Brownian motion. Hara and Slade
(2000a,b) studied the convergence of rescaled percolation clusters
to integrated super-Brownian excursions. Those results provide
interesting connections between superprocesses and interacting
particle systems.

A class of superprocesses with dependent spatial motion (SDSM)
over the real line $\dfR$ were introduced and constructed in Wang
(1997, 1998). The construction was then generalized in Dawson {\it
et al.} (2001). In the model, the spatial motion is defined by a
system of differential equations driven by a family of independent
Brownian motions, the individual noises, and a time-space white
noise, the common noise. If the coefficient of the individual
noises are uniformly bounded away from zero, the SDSM is
absolutely continuous and its density satisfies a stochastic
differential equation (SPDE) that generalizes the Konno-Shiga
equation satisfied by super Brownian motion over $\dfR$; see
Dawson {\it et al.} (2000) and Konno and Shiga (1988). When the
individual noises vanish, the SDSM is purely atomic; see Wang
(1997, 2002). A construction of the purely atomic SDSM in terms of
one-dimensional excursions was given in Dawson and Li (2003),
where an immigration diffusion process was also constructed as the
pathwise unique solution of a stochastic equation. An SPDE for the
purely atomic SDSM was derived recently in Dawson {\it et al.}
(2003). It was proved in Dawson {\it et al.} (2001) that a
suitably rescaled absolutely continuous SDSM converges to the
usual super Brownian motion. This describes another situation
where the super Brownian motion arises universally. For the purely
atomic SDSM, it was mentioned in the introduction of Dawson {\it
et al.} (2001) that the same rescaled limit would lead to a
superprocess with coalescing spatial motion (SCSM), a continuous
state version of the coalescing-branching particle system. This
seems to be a new phenomenon in scaling limits of interacting
particle systems and superprocesses. The statement was not proved
in Dawson {\it et al.} (2001) since the construction and
characterization of the SCSM remained open at that time.

The main purpose of this paper is to give a proof of the
observation of Dawson {\it et al.} (2001). As a preliminary, we
give in Section 2 some characterizations for a coalescing Brownian
flow in terms of martingale problems and show that the flow is
actually the scaling limit of the interacting Brownian flow that
serves as the carrier of the purely atomic SDSM in the excursion
representation given in Dawson and Li (2003).

In Section 3, we construct the SCSM from the coalescing Brownian
flow and one-dimensional excursions following the idea of Dawson
and Li (2003). It has been known for a long time that a
superprocess without spatial motion reduces to a Poisson system of
point masses that evolve according to Feller branching diffusions
without interaction; see Shiga (1990). The SCSM adds a coalescing
Brownian flow which carries the point masses. Any masses join
together when their carriers coalesce.

In Section 4, we derive the scaling limit theorem of the SDSM from
that of the interacting Brownian flow and the excursion
representations. This result shows that excursion representations
play important roles not only in the construction of the
superprocesses but also in the study of some of their properties.

\section{Interacting Brownian flows}

\setcounter{equation}{0}

An $m$-dimensional continuous process $\{(y_1(t), \cdots, y_m(t)):
t\ge0\}$ is called an {\it $m$-system of coalescing Brownian
motions} ($m$-SCBM) with speed $\rho>0$ if each $\{y_i(t): t\ge
0\}$ is a Brownian motion with speed $\rho>0$ and, for $i\neq j$,
$\{|y_i(t) - y_j(t)|: t\ge0\}$ is a Brownian motion with speed
$2\rho$ stopped at the origin. Clearly, $\{(y_1(t), \cdots,
y_m(t)): t\ge0\}$ is an $m$-SCBM with speed $\rho>0$ if and only
if
 \beqlb\label{2.1}
\<y_i,y_j\>(t) = \rho \cdot (t-t\land \tau_{ij}),
\qquad 1\le i, j\le m,
 \eeqlb
where $\tau_{ij} = \inf\{t\ge0: y_i(t) = y_j(t)\}$.

To give a martingale characterization of the $m$-SCBM, we need to
choose a convenient core of its generator. For any permutation
$(i_1,i_2,\cdots,i_m)$ of $(1,2,\cdots,m)$ let
 \beqlb\label{2.2}
\dfR^m_{i_1i_2\cdots i_m} = \{(x_1,x_2,\cdots,x_m)\in \dfR^m:
x_{i_1} <x_{i_2} <\cdots <x_{i_m}\}.
 \eeqlb
Let ${\mcr D}_0^{(1)} = C^2(\dfR)$ and for $m\ge 2$ let ${\mcr
D}^{(m)}_0$ be the set of functions $f\in C(\dfR^m)$ such that

\begin{itemize}

\item[(2.A)] $f$ is twice continuously differentiable in each
$\dfR^m_{i_1i_2 \cdots i_m}$ with bounded partial derivatives up
to the second degree;

\item[(2.B)] all partial derivatives of $f$ up to the second
degree can be extended to the closure of each $\dfR^m_{i_1i_2
\cdots i_m}$ as uniformly continuous functions with
 \beqlb\label{2.3}
\frac{\partial^2f}{\partial x_i\partial x_j} (x_1,\cdots,x_m) =
\frac{\partial^2f}{\partial x_j\partial x_i} (x_1,\cdots,x_m) = 0
 \eeqlb
for $1\le i<j\le m$ and $x_i=x_j$. (We simply write
$\frac{\partial^2f}{\partial x_i\partial x_j}$ for the continuous
extension of the derivative.)

\end{itemize}

For any $1\le i< j\le m$, define the operator $p^{(m)}_{ij}:
C(\dfR^m) \to C(\dfR^{m-1})$ by
 \beqnn
p_{ij}f(x_1,\cdots,x_{m-1}) =
f(x_1,\cdots,x_{m-1},\cdots,x_{m-1},\cdots,x_{m-2}),
 \eeqnn
where $x_{m-1}$ occurs at the places of the $i$th and the $j$th
variables on the right hand side. Let ${\mcr D}^{(m)}$ be the
totality of functions $f\in C(\dfR^m)$ such that
$p^{(m-k)}_{i_kj_k} \cdots p^{(m)}_{i_0j_0} f\in {\mcr D}
^{(m-k-1)}_0$ for all $1\le i_l<j_l\le m-l$ and $1\le l\le m-1$.
For $f\in {\mcr D}^{(m)}$, let
 \beqlb\label{2.4}
G_0^{(m)} f(x_1,\cdots,x_m) = \frac{1}{2}\rho \itDelta_m
f(x_1,\cdots,x_m), \quad (x_1,\cdots,x_m)\in \dfR^m,
 \eeqlb
where $\itDelta_m$ denotes the $m$-dimensional Laplace operator.

A continuous process $\{(x_1(t),\cdots,x_m(t)): t\ge0\}$ is called
a solution of the {\it $(G_0^{(m)}, {\mcr D}^{(m)})$-martingale
problem} if
 \beqnn
f(x_1(t),\cdots,x_m(t)) - f(x_1(0),\cdots,x_m(0)) - \int_0^t
G_0^{(m)} f(x_1(s),\cdots,x_m(s))ds
 \eeqnn
is a martingale for each $f\in {\mcr D}^{(m)}$.

\btheorem\label{t2.1} The distribution $\bfP_{(b_1,b_2)}$ on
$C([0,\infty), \dfR^2)$ of the $2$-SCBM with speed $\rho>0$ and
initial state $(b_1,b_2)$ is the unique probability measure on
$C([0,\infty), \dfR^2)$ such that $\bfP_{(b_1,b_2)} \{(w_1(0),$ $
w_2(0)) = (b_1,b_2)\} =1$ and $\{(w_1(t),w_2(t)): t\ge 0\}$ under
$\bfP_{(b_1,b_2)}$ solves the $(G_0^{(2)}, {\mcr
D}^{(2)})$-martingale problem. \etheorem

{\it Proof.} We first show that the $2$-SCBM solves the
$(G_0^{(2)}, {\mcr D}^{(2)})$-martingale problem. Let $f\in {\mcr
D}^{(2)}$. If $b_1=b_2$, then $\bfP_{(b_1,b_2)} \{w_1(t)=w_2(t)$
for all $t\ge0\}=1$. By It\^o's formula we get
 \beqlb\label{2.5}
& & f(w_1(t),w_2(t)) - f(w_1(0),w_2(0)) \nonumber \\
&=&
\sum_{i=1}^2\int_0^t f^{\prime}_i(w_1(s),w_2(s))dw_i(s)
+ \frac{1}{2}\rho\sum_{i,j=1}^2\int_0^t
f^{\prime\prime}_{ij} (w_1(s),w_2(s))d\<w_i,w_j\>(s) \nonumber \\
&=&
\sum_{i=1}^2\int_0^t f^{\prime}_i(w_1(s),w_2(s))dw_i(s)
+ \frac{1}{2}\rho\sum_{i=1}^2\int_0^t f^{\prime\prime}_{ii}
(w_1(s),w_2(s)) ds,
 \eeqlb
where we have used the assumption $f^{\prime\prime}_{12}
(x,x) = f^{\prime\prime}_{21}(x,x)=0$ for the last equality.
If $b_1\neq b_2$, we have
 \beqlb\label{2.6}
& &f(w_1(t\land\tau),w_2(t\land\tau)) - f(w_1(0),w_2(0))
\nonumber \\
&=&
\sum_{i=1}^2\int_0^{t\land\tau} f^{\prime}_i
(w_1(s),w_2(s))dw_i(s)
+ \frac{1}{2}\rho\sum_{i=1}^2\int_0^{t\land\tau}
f^{\prime\prime}_{ii} (w_1(s),w_2(s)) ds,
 \eeqlb
where $\tau = \inf\{t\ge0: w_1(t) = w_2(t)\}$. Summing up
(\ref{2.5}) and (\ref{2.6}) we see that $\{(w_1(t),w_2(t)): t\ge
0\}$ under $\bfP_{(b_1,b_2)}$ is a solution of the $(G_0^{(2)},
{\mcr D}^{(2)})$-martingale problem. Conversely, suppose that
$\bfP_{(b_1,b_2)}$ is a probability measure on $C([0,\infty),
\dfR^2)$ such that $\bfP_{(b_1,b_2)} \{(w_1(0),w_2(0)) =
(b_1,b_2)\} =1$ and $\{(w_1(t),w_2(t)): t\ge 0\}$ under
$\bfP_{(b_1,b_2)}$ solves the $(G_0^{(2)}, {\mcr
D}^{(2)})$-martingale problem. For $f\in C^2(\dfR)$ we apply the
martingale problem to the function $(y_1,y_2) \mapsto f(y_1)$ to
see that
 \beqnn
f(w_1(t)) - f(y_1)
=
\mbox{mart.} + \frac{1}{2}\rho\int_0^t
f^{\prime\prime} (w_1(s)) ds.
 \eeqnn
Therefore, $\{w_1(t): t\ge 0\}$ under $\bfP_{(b_1,b_2)}$ is a
Brownian motion with speed $\rho$. Similarly, $\{w_2(t): t\ge 0\}$
under $\bfP_{(y_1,y_2)}$ is also a Brownian motion with speed
$\rho$. On the other hand, for $f\in C^2_0([0,\infty))$ satisfying
$f^{\prime\prime}(0)=0$ we find by applying the martingale problem
to the function $(y_1,y_2) \mapsto f(|y_1-y_2|)$ that
 \beqlb\label{2.7}
f(|w_1(t)-w_2(t)|) - f(|y_1-y_2|)
=
\mbox{mart.} + \rho\int_0^t
f^{\prime\prime} (|w_1(s)-w_2(s)|) ds.
 \eeqlb
For $n\ge 1$ let $f_n \in C^2_0([0,\infty))$ be such that $f_n(x)
= x$ for $0 \le x\le n$. Applying (\ref{2.7}) to the sequence
$\{f_n\}$ we see that $\{|w_1(t) - w_2(t)|: t\ge 0\}$ under
$\bfP_{(b_1,b_2)}$ is a non-negative local martingale so it must
be absorbed at zero. By It\^o's formula,
 \beqnn
|w_1(t)-w_2(t)|^2 - |y_1-y_2|^2 = \mbox{local mart.} + \int_0^t
d\<|w_1-w_2|\>(s).
 \eeqnn
For $n\ge 1$ let $h_n \in C^2_0([0,\infty))$ be such that
$h_n^{\prime\prime} (0) = 0$ and $h_n(x) = x^2$ for $n^{-1} \le
x\le n$. Since $\{|w_1(t) - w_2(t)|: t\ge 0\}$ is absorbed at
zero, applying (\ref{2.7}) to $\{h_n\}$ we see that
 \beqnn
|w_1(t)-w_2(t)|^2 - |y_1-y_2|^2 = \mbox{local mart.} + 2\rho
(t\land \tau),
 \eeqnn
where $\tau = \inf \{t\ge0: w_1(t) = w_2(t)\}$. It follows that
$\<|w_1-w_2|\>(t) = 2\rho \cdot (t\land \tau)$ and hence
$\bfP_{(b_1,b_2)}$ is the distribution of the $2$-SCBM. \qed

\btheorem\label{t2.2} The distribution $\bfP_{(b_1,\cdots,b_m)}$
on $C([0,\infty), \dfR^m)$ of the $m$-SCBM with speed $\rho>0$ and
initial state $(b_1,\cdots,b_m)$ is the unique probability measure
on $C([0,\infty), \dfR^m)$ such that $\bfP_{(b_1,\cdots,b_m)}
\{(w_1(0),\cdots,w_m(0)) = (b_1,\cdots,b_m)\} =1$ and
$\{(w_1(t),\cdots,w_m(t)): t\ge 0\}$ under
$\bfP_{(b_1,\cdots,b_m)}$ solves the $(G_0^{(m)}, {\mcr
D}^{(m)})$-martingale problem. \etheorem

{\it Proof.} By considering the $m$-SCBM piece by piece between
the coalescing times as in the proof of the last theorem, one can
show that it is indeed a solution of the $(G_0^{(m)}, {\mcr
D}^{(m)})$-martingale problem. To see the uniqueness, observe that
for any $1\le i<j\le m$ and $f\in {\mcr D}^{(2)}$, the function
$(y_1,\cdots,y_m) \mapsto f(y_i,y_j)$ belongs to ${\mcr D}^{(m)}$.
It follows that if $\{(w_1(t),\cdots,w_m(t)): t\ge 0\}$ under
$\bfP_{(b_1,\cdots,b_m)}$ is a solution of the $(G_0^{(m)}, {\mcr
D}^{(m)})$-martingale problem, then the pair $\{(w_i(t),w_j(t)):
t\ge 0\}$ is a solution of the $(G_0^{(2)}, {\mcr
D}^{(2)})$-martingale problem. By Theorem~\ref{t2.1},
$\{(w_i(t),w_j(t)): t\ge 0\}$ under $\bfP_{(b_1,\cdots,b_m)}$ is a
$2$-SCSM and hence $\<w_i,w_j\>(t) = \rho \cdot (t-t\land
\tau_{ij})$, where $\tau_{ij} = \inf\{t\ge0: w_i(t) = w_j(t)\}$.
Then $\{(w_1(t),\cdots,w_m(t)): t\ge 0\}$ is an $m$-SCBM. \qed

We may embed the $m$-SCBM into an inhomogeneous Markov process
with state space $W := C([0,\infty)$, $\dfR)$. To this end, let
$W^{\dfR}$ denote the totality of $W$-valued paths $\{w(a,\cdot):
a\in \dfR\}$, which contains all possible paths of the Markov
process to be defined. For any $\{b_1, \cdots, b_m\}\subset \dfR$,
let $F_{b_1,\cdots,b_m}$ denote the distribution on $W^m$ of the
$m$-SCBM $\{(y_1(t), \cdots, y_m(t)): t\ge0\}$ with speed $\rho>0$
and initial state $(b_1, \cdots, b_m)$. It is easy to see that
$\{F_{b_1,\cdots,b_m}: b_1, \cdots, b_m \in \dfR\}$ is a
consistent family. By Kolmogorov's extension theorem, there is a
unique probability measure $\bfP^{cb}$ on $W^{\dfR}$ which has
finite dimensional distributions $\{F_{b_1,\cdots,b_m}: b_1,
\cdots, b_m \in \dfR\}$. A two parameter process $\{y(a,t): a\in
\dfR, t\ge0\}$ is called a {\it coalescing Brownian flow} if the
path-valued process $\{y(a,\cdot): a\in \dfR\}$ has distribution
$\bfP^{cb}$ on $W^{\dfR}$. Indeed, $\{y(a,\cdot): a\in \dfR\}$ is
an inhomogeneous Markov process. For $a\in \dfR$ let $W_a$ denote
the set of paths $w\in W$ with $w(0)=a$. For any $a\in \dfR$ let
$\{B_a(t): t\ge0\}$ be a Brownian motion with speed $\rho>0$ and
initial state $B_a(0)=a$ and let $Q_a(\cdot)$ denote the
distribution of $\{B_a(t): t\ge0\}$ on $W_a$. For $a\le b\in \dfR$
and $w_a\in W_a$ let $Q_{a,b} (w_a,\cdot)$ denote the distribution
on $W_b$ of the process $\{\xi_b(t): t\ge0\}$ defined by
 \beqnn
\xi_b(t) =\left\{\begin{array}{ll}
B_b(t)      \quad&\mbox{ if $t\le\tau_{ab}$,}  \\
w_a(t)           &\mbox{ if $t>\tau_{ab}$,}
\end{array}\right.
 \eeqnn
where $\tau_{ab}= \inf\{t\ge0: B_b(t)=w_a(t)\}$. Then $(Q_{a,b})
_{a\le b}$ is a Markov transition semigroup with state spaces
$\{W_a: a\in \dfR\}$ and $(Q_a)_{a\in\dfR}$ is an entrance law for
$(Q_{a,b})_{a\le b}$. It is not hard to see that a coalescing
Brownian flow $\{y(a,\cdot): a\in \dfR\}$ is a Markov process with
transition semigroup $(Q_{a,b})_{a\le b}$ and one-dimensional
distributions $(Q_a)_{a\in\dfR}$. See Dynkin (1978, p.724) for
discussions of inhomogeneous Markov processes determined by
entrance laws. A more general coalescing Brownian flow is defined
and studied in Harris (1984), where interaction is allowed between
the particles before they coalesce. We refer the reader to Evans
and Pitman (1998) and the references therein for some recent work
on related models.

Now we consider an interacting Brownian flow driven by a
time-space white noise. Let $h\in C(\dfR)$ be square-integrable
and continuously differentiable with square-integrable derivative
$h^\prime$. Suppose we are given on some standard probability
space $(\itOmega,{\mcr F},\bfP)$ a time-space white noise
$W(ds,dy)$ on $[0,\infty)\times\dfR$ based on the Lebesgue
measure; see, e.g., Walsh (1986). By Dawson {\it et al.} (2001)
and Wang (1997, 1998), for each $a\in\dfR$ the equation
 \beqlb\label{2.8}
x(t) = a + \int_0^t\int_{\dfR} h(y-x(s))W(ds,dy), \quad t\ge 0,
 \eeqlb
has a unique solution $\{x(a,t): t\ge0\}$. We call $\{x(a,t):
t\ge0; a\in \dfR\}$ an {\it interacting Brownian flow} driven by
the time-space white noise. It is not hard to check that for any
$(a_1,\cdots,a_m) \in \dfR^m$, the solutions $\{(x(a_1,t), \cdots,
x(a_m,t)): t\ge0\}$ of (\ref{2.8}) constitute an $m$-dimensional
diffusion process generated by the differential operator
 \beqlb\label{2.9}
G^{(m)} := \frac{1}{2}\sum_{i,j=1}^m
\rho(x_i-x_j)\frac{\partial^{2}}{\partial x_i \partial x_j},
 \eeqlb
where
 \beqlb\label{2.10}
\rho(x) = \int_{\dfR}h(y-x)h(y) dy, \quad x\in\dfR.
 \eeqlb
In particular, each $\{x(a_i,t): t\ge0\}$ is a one-dimensional
Brownian motion with quadratic variation process $\rho(0)t$, so we
call $\{(x(a_1,t),\cdots,x(a_n,t)): t\ge0\}$ an {\it $m$-system of
interacting Brownian motions} ($m$-SIBM).

Given $\theta>0$ and $f\in C(\dfR)$, let $f_\theta(x) = f(\theta
x)$. Replacing $h(\cdot)$ in (\ref{2.8}) by $\sqrt {\theta}
h_\theta (\cdot)$ we obtain the function $\rho_\theta(\cdot)$, so
the latter can also serve as the interaction parameter of an
interacting Brownian motion. The following theorem shows that the
coalescing Brownian flow arises in some sense as the scaling limit
of the interacting Brownian flow driven by the time-space white
noise.

\btheorem\label{t2.3} Suppose that $\rho(x)\to 0$ as
$|x|\to\infty$. For each $\theta \ge1$, let $\{(x^{\theta}_1(t),
\cdots, x^{\theta}_m(t)): t\ge0\}$ be an $m$-SIBM with interaction
parameter $\rho_\theta(\cdot)$ and initial state $(a^{\theta}_1,
\cdots, a^{\theta}_m)$. If $a^{\theta}_i \to b_i$ as $\theta \to
\infty$, then the law of $\{(x^{\theta}_1(t), \cdots,
x^{\theta}_m(t)): t\ge0\}$ on $C([0,\infty),\dfR^m)$ converges to
that of the $m$-SCBM with speed $\rho(0)$ starting from $(b_1,
\cdots, b_m)$. \etheorem

{\it Proof.} The result could be proved using Theorem~\ref{t2.2}.
The following proof directly based on the definition of the SIBM
seems more readable. Since each $x_i^{\theta}(t)$ is a Brownian
motion with speed $\rho(0)$, we get by Doob's martingale
inequality that
 \beqnn
\bfP\bigg\{\sup_{0\le t\le T} |x_i^{\theta}(t)| > \eta\bigg\} \le
\bfP\{x_i^{\theta}(T)^2\}/\eta^2 = 2(|a_i^{\theta}|^2 +
\rho(0)T)/\eta^2,
 \eeqnn
where we also use $\bfP$ to denote the expectation; see, e.g.,
Ikeda and Watanabe (1989, p.34). Then for each $\varepsilon > 0$
there is a compact set $K_\varepsilon \subset \dfR^m$ such that
 \beqlb\label{2.11}
\sup_{\theta\ge 1} \bfP\{(x^{\theta}_1(t), \cdots,
x^{\theta}_m(t)) \in K_\varepsilon \mbox{ for } 0\le t\le T \} \le
\varepsilon,
 \eeqlb
that is, the family $\{(x^{\theta}_1(\cdot), \cdots, x^{\theta}_m
(\cdot)): \theta\ge 1\}$ satisfies the compact containment
condition of Ethier and Kurtz (1986, p.142). Let
 \beqnn
G_\theta^{(m)} := \frac{1}{2}\sum_{i,j=1}^m \rho_\theta(x_i-x_j)
\frac{\partial^{2}}{\partial x_i \partial x_j}.
 \eeqnn
Let $C_\kappa(\dfR^m)$ denote the set of continuous functions on
$\dfR$ with compact supports and let $H = C_\kappa(\dfR^m) \cap
C^2(\dfR^m)$. Then for $f\in H$,
 \beqlb\label{2.12}
f(x_1^{\theta}(t),\cdots,x_m^{\theta}(t)) -
f(a_1^{\theta},\cdots,a_m^{\theta}) - \int_0^t G_\theta^{(m)}
f(x_1^{\theta}(s),\cdots,x_m^{\theta}(s)) ds
 \eeqlb
is a martingale. Observe that $\sup_{\theta\ge1} \|G_\theta^{(m)}
f\| < \infty$, so for each $T>0$ we have
 \beqnn
\sup_{\theta\ge1}\bfE\bigg[\bigg(\int_0^T|G_\theta^{(m)}
f(x_1^{\theta}(s),\cdots,x_m^{\theta}(s))|^2ds\bigg)^{1/2}\bigg] <
\infty.
 \eeqnn
By Ethier and Kurtz (1986, p.145), $\{f(x^{\theta}_1(\cdot),
\cdots, x^{\theta}_m(\cdot)): \theta\ge 1\}$ is a tight family in
$C([0,\infty), \dfR)$, which is a closed subspace of
$D([0,\infty), \dfR)$. Since $H$ is dense in $C_\kappa(\dfR^m)$ in
the topology of uniform convergence on compact sets, by Ethier and
Kurtz (1986, p.142), $\{(x^{\theta}_1(\cdot), \cdots,
x^{\theta}_m(\cdot)): \theta\ge 1\}$ is tight in $C([0,\infty),
\dfR^m)$. Let $\bfP_0$ be the limit distribution on $C([0,\infty),
\dfR^m)$ of any convergent subsequence $(x^{\theta_k}_1(t),
\cdots, x^{\theta_k}_m(t))$ with $\theta_k \to \infty$. Since each
$x_i^{\theta}(t)$ is a Brownian motion with speed $\rho(0)$, so is
$w_i(t)$ under $\bfP_0$. As in Wang (1998, p.756), one may see
that $\{x_i^{\theta}(t) - x_j^{\theta}(t): t\ge0\}$ is a diffusion
process for which the origin is an unaccessible trap. It follows
that $\bfP\{x_i^{\theta} (t) = x_j^{\theta} (t)$ for all $t\ge0\}
= 1$ if $a_i^{\theta} = a_j^{\theta}$ and $\bfP\{x_i^{\theta}(t)
\neq x_j^{\theta}(t)$ for all $t\ge0\} = 1$ if $a_i^{\theta} \neq
a_j^{\theta}$. In view of (\ref{2.12}), for any $f\in
C^2_0([0,\infty))$ with $f ^{\prime\prime} (0) =0$,
 \beqlb\label{2.13}
&& f(|x_j^{\theta}(t)-x_i^{\theta}(t)|)
- f(|a_j^{\theta}-a_i^{\theta}|)  \nonumber \\
&&\hskip2cm  -\int_0^t [\rho(0) - \rho_\theta(x_j^{\theta}(s)
-x_i^{\theta}(s))]
f^{\prime\prime}(|x_j^{\theta}(s) - x_i^{\theta}(s)|) ds
 \eeqlb
is a martingale. Since $f ^{\prime\prime} (0) =0$ and $\rho(x)\to
0$ as $|x|\to\infty$, we have $[\rho(0) - \rho_\theta(\cdot)]
f^{\prime\prime}(|\cdot|) \to \rho(0)f^{\prime\prime}(|\cdot|)$
uniformly as $\theta \to \infty$. Letting $\theta \to \infty$ in
(\ref{2.13}) along $\{\theta_k\}$ we see
 \beqnn
f(|w_j(t)-w_i(t)|) - f(|b_j-b_i|) - \int_0^t
\rho(0)f^{\prime\prime}(|w_j(s) - w_i(s)|) ds
 \eeqnn
under $\bfP_0$ is a martingale. As in the proof of
Theorem~\ref{t2.1}, $\{|w_j(t) - w_i(t)|: t\ge0\}$ under $\bfP_0$
must be a non-negative local martingale having quadratic variation
process $2\rho(0) (t\land \tau_{ij})$ with $\tau_{ij} = \inf
\{t\ge0: w_i(t) = w_j(t)\}$. Thus $\bfP_0$ is the law of the
$m$-SCBM Brownian motion starting from $(b_1,\cdots,b_m)$ with
speed $\rho(0)$. \qed

\section{Superprocesses with coalescing spatial motion}

\setcounter{equation}{0}

In this section, we give some constructions for the SCSM. Let
$\rho >0$ be a constant. Suppose that $\sigma \in C(\dfR)^+$ and
there is a constant $\epsilon >0$ such that $\inf_x\sigma (x) \ge
\epsilon$. The formal generator of the SCSM is given by
 \beqlb\label{3.1}
\L F(\mu) &=& \frac{1}{2}\int_{\dfR} \sigma\frac{\delta^2 F(\mu)}
{\delta\mu(x)^2}\mu(dx)
+\,\frac{1}{2}\rho\int_{\dfR}\frac{d^2}{dx^2}
\frac{\delta F(\mu)}{\delta\mu(x)}\mu(dx)  \nonumber  \\
& &
+\,\frac{1}{2}\int_{\itDelta}\frac{d^2}{dxdy}
\frac{\delta^2 F(\mu)}{\delta\mu(x)\delta\mu(y)}\mu(dx)\mu(dy),
 \eeqlb
where $\itDelta = \{(x,x): x\in \dfR\}$. Note that the first two
terms on the right hand side simply give the generator of a usual
super Brownian motion, where the first term describes the
branching and the second term gives the spatial motion. The last
term shows that interactions in the spatial motion only occur
between `particles' located at the same positions. Those
descriptions are justified by the constructions to be given.

We first consider a purely atomic initial state with a finite
number of atoms. In the sequel, a martingale diffusion $\{\xi(t):
t\ge0\}$ is called a {\it standard Feller branching diffusion} if
it has quadratic variation $\xi(t)dt$. Let $\{(y_1(t), \cdots,
y_n(t)): t\ge0\}$ be an $n$-SCBM with speed $\rho$ and initial
state $(b_1, \cdots, b_n)\in \dfR^n$. Let $\{(\xi_1(t), \cdots,
\xi_n(t)): t\ge0\}$ be a system of independent standard Feller
branching diffusions with initial state $(\xi_1, \cdots, \xi_n)\in
\dfR_+^n$. We assume that $\{(y_1(t),\cdots,y_n(t)): t\ge0\}$ and
$\{(\xi_1(t), \cdots, \xi_n(t)): t\ge0\}$ are defined on a
standard complete probability space $(\itOmega, {\mcr F}, \bfP)$
and are independent of each other. Set
 \beqlb\label{3.2}
\psi^\sigma_i(t) = \int_0^t \sigma(y_i(s)) ds
 \eeqlb
and $\xi_i^\sigma(t) = \xi_i(\psi^\sigma_i(t))$. Then
 \beqlb\label{3.3}
X_t = \sum_{i=1}^n \xi^\sigma_i(t)\delta_{y_i(t)}, \quad t\ge0,
 \eeqlb
defines a continuous $M(\dfR)$-valued process. Intuitively, this
process consists of $n$ particles carried by the $n$-SCBM
$\{(y_1(t),\cdots,y_n(t)): t\ge0\}$. The mass of the $i$th
particle is given by $\{\xi^\sigma_i(t): t \ge 0\}$, which is
obtained from a standard Feller branching diffusion by a time
change depending on the position of the $i$th carrier. Thus we
have here a spatially dependent branching mechanism.

Let ${\mcr G}_t$ be the $\sigma$-algebra generated by the family
of $\bfP$-null sets in ${\mcr F}$ and the family of random
variables
 \beqlb\label{3.4}
\{(y_1(s), \cdots, y_n(s),\xi_1^\sigma(s), \cdots, \xi_n^\sigma(s)):
0\le s \le t\}.
 \eeqlb
Then we have

\btheorem\label{t3.1} The process $\{X_t: t\ge0\}$ defined by
(\ref{3.3}) is a diffusion process relative to the filtration
$({\mcr G}_t)_{t\ge0}$ with state space $M_a(\dfR)$, the set of
purely atomic measures on $\dfR$. \etheorem

{\it Proof.} Let $\mu= \sum_{i=1}^n \xi_i\delta_{a_i}$. By
symmetry, the distribution $Q_t(\mu,\cdot)$ of $X_t$ on
$M_a(\dfR)$ only depends on $t\ge0$ and $\mu$. Clearly, under
$\bfP\{\cdot\, | {\mcr G}_r\}$ the process $\{(x_1(r+t), \cdots,
x_n(r+t)): t\ge 0\}$ is an $n$-SCBM and
$\{(\xi_1(\psi^\sigma_i(r)+t), \cdots, \xi_n(\psi^\sigma_i(r)+t)):
t\ge 0\}$ is a system of independent standard Feller branching
diffusions. Moreover, the two systems are conditionally
independent of each other. Then $X_{r+t}$ under $\bfP\{\cdot\, |
{\mcr G}_r\}$ has distribution $Q_t (X_r,\cdot)$. The Feller
property of the $Q_t(\mu,\cdot)$ follows from those of
$(x_1(t),\cdots,x_n(t))$ and $(\xi_1(t), \cdots, \xi_n(t))$. Then
the strong Markov property holds by the continuity of $\{X_t:
t\ge0\}$. \qed

\btheorem\label{t3.2} If $\{X_t: t\ge0\}$ is given by (\ref{3.3}),
then for each $\phi\in C^2(\dfR)$,
 \beqlb\label{3.5}
M_t(\phi) = \<\phi,X_t\> - \<\phi,X_0\>
- \frac{1}{\,2\,}\rho \int_0^t \<\phi^{\prime\prime},X_s\> ds,
\quad t\ge0,
 \eeqlb
is a continuous martingale relative to $({\mcr G}_t)_{t\ge0}$ with
quadratic variation process
 \beqlb\label{3.6}
\<M(\phi)\>_t = \int_0^t\<\sigma\phi^2,X_s\> ds + \int_0^t
ds\int_{\itDelta}\phi^\prime(x)\phi^\prime(y) X_s(dx)X_s(dy),
 \eeqlb
where $\itDelta = \{(x,x): x\in \dfR\}$. \etheorem

\noindent{\it Proof.} As in the proof of Dawson and Li (2003,
Theorem~3.3), $\{(\xi_i^\sigma(t): t\ge 0\}$ is a continuous
martingale with quadratic variation $\sigma(y_i(t))dt$ and
$\<\xi_i^\sigma, \xi_j^\sigma\>(t) \equiv 0$ if $i\neq j$. By
It\^o's formula,
 \beqlb\label{3.7}
\xi_i^\sigma(t)\phi(y_i(t))
&=& \xi_i^\sigma(0)\phi(y_i(0))
+ \int_0^t\phi(y_i(s)) d\xi_i^\sigma(s)
+ \int_0^t\xi_i^\sigma(s)\phi^\prime(y_i(s)) dy_i(s)  \nonumber   \\
& & + \frac{1}{\,2\,}\rho \int_0^t\xi_i^\sigma(s)
\phi^{\prime\prime} (y_i(s)) ds.
 \eeqlb
Taking the summation we get
 \beqnn
\<\phi,X_t\> - \<\phi,X_0\>
= M_t(\phi)
+ \frac{1}{\,2\,}\rho \int_0^t \<\phi^{\prime\prime},X_s\> ds,
\quad t\ge0,
 \eeqnn
where
 \beqnn
M_t(\phi)
:=
\sum_{i=1}^n \int_0^t\phi(y_i(s))d\xi_i^\sigma(s)
+ \sum_{i=1}^n \int_0^t \xi_i^\sigma(s)\phi^\prime(y_i(s))dy_i(s),
 \eeqnn
is a continuous martingale relative to $({\mcr G}_t)_{t\ge0}$ with
quadratic variation process
 \beqnn
\<M(\phi)\>_t
&=&
\sum_{i=1}^n \int_0^t \sigma(y_i(s))\xi_i^\sigma(s)\phi(y_i(s))^2 ds
+ \sum_{i,j=1}^n \int_{\tau_{ij}}^t \xi_i^\sigma(s)\xi_j^\sigma(s)
\phi^\prime(y_i(s))\phi^\prime(y_j(s))ds  \\
&=&
\int_0^t\<\sigma\phi^2,X_s\> ds
+ \int_0^t ds\int_{\itDelta}\phi^\prime(x)\phi^\prime(y)
X_s^2(dx,dy),
 \eeqnn
where $\tau_{ij} = \inf\{s\ge0: y_i(s) = y_j(s)\}$. This gives the
desired result. \qed

We can give another martingale characterization of the process as
follows. Let ${\mcr D}(\L)$ be the set of all functions of the
form $F_{m,f}(\mu) = \<f,\mu^m\>$ with $f\in {\mcr D}^{(m)}$.
Observe that
 \beqlb\label{3.8}
\L F_{m,f}(\mu) = F_{m,G_0^{(m)}f}(\mu) + \frac{1}{2} \sum_{i,j=1,
i\neq j}^m F_{m-1,\itPhi_{ij}f}(\mu),
 \eeqlb
where $G_0^{(m)}$ is the generator of the $m$-SCBM with speed
$\rho$ and $\itPhi_{ij}$ denotes the operator from $C(\dfR^m)$ to
$C(\dfR^{m-1})$ defined by
 \beqlb\label{3.9}
\itPhi_{ij}f(x_1,\cdots,x_{m-1})
=
\sigma(x_{m-1})f(x_1,\cdots,x_{m-1},\cdots,x_{m-1},\cdots,x_{m-2}),
 \eeqlb
where $x_{m-1}$ takes the places of the $i$th and the $j$th variables of
$f$ on the right hand side.

\btheorem\label{t3.3} Let $\{X_t: t\ge0\}$ be defined by
(\ref{3.3}). Then $\bfE \{\<1,X_t\>^m\}$ is locally bounded in
$t\ge0$ for each $m\ge1$ and $\{X_t: t\ge0\}$ solves the $({\mcr
L}, {\mcr D}(\L))$-martingale problem, that is, for each $F_{m,f}
\in {\mcr D}(\L)$,
 \beqlb\label{3.10}
F_{m,f}(X_t) - F_{m,f}(X_0) - \int_0^t {\mcr L}F_{m,f}(X_s)ds
 \eeqlb
is a martingale. \etheorem

{\it Proof.} Based on Theorem~\ref{t3.2}, it is not hard to show
that $\bfE \{\<1,X_t\>^m\}$ is locally bounded in $t\ge0$ for each
$m\ge1$. Since $\{(\xi_i^\sigma(t): t\ge 0\}$ is a continuous
martingale with quadratic variation $\sigma(y_i(t))dt$ and
$\<\xi_i^\sigma, \xi_j^\sigma\>(t) \equiv 0$ if $i\neq j$, for
$m\ge1$ and $f\in {\mcr D}^{(m)}$ we have
 \beqnn
\<f,X_t^m\>
&=& \sum_{i_1,\cdots,i_m=1}^n \xi_{i_1}^\sigma(t) \cdots
\xi_{i_m}^\sigma(t) f(y_{i_1}(t),\cdots,y_{i_m}(t))   \\
&=& \sum_{i_1,\cdots,i_m=1}^n \xi_{i_1}^\sigma(0) \cdots
\xi_{i_m}^\sigma(0) f(y_{i_1}(0),\cdots,y_{i_m}(0))
+ \mbox{mart.} \\
&& + \,\frac{1}{\,2\,}\rho\sum_{i_1,\cdots,i_m=1}^n\sum_{j=1}^m
\int_0^t\xi_{i_1}^\sigma(s) \cdots \xi_{i_m}^\sigma(s)
f^{\prime\prime}_{jj}(y_{i_1}(s),\cdots,y_{i_m}(s))ds \\
&& +
\,\frac{1}{\,2\,}\sum_{\alpha,\beta=1}^m\sum_{\{\mbox{cond.}\}}
\int_0^t\sigma(y_{i_\alpha}(s))\xi_{i_1}^\sigma(s) \cdots
\xi_{i_m}^\sigma(s)\xi_{i_\alpha}^\sigma(s)^{-1}
f(y_{i_1}(s),\cdots,y_{i_m}(s))ds \\
&=& \<f,X_0^m\> + \mbox{mart.}
+ \,\frac{1}{\,2\,}\rho \int_0^t \< \itDelta f,X_s^m\>ds
+ \,\frac{1}{\,2\,}\sum_{\alpha,\beta=1}^m
\int_0^t \< \itPhi_{\alpha\beta} f,X_s^{m-1}\>ds,
 \eeqnn
where $\{\mbox{cond.}\} = \{$ for all $1\le i_1,\cdots,i_m\le n$
with $i_\alpha=i_\beta \}$ and we used the fact
$f^{\prime\prime}_{ij} (x_1,\cdots,x_m) = 0$ for $x_i = x_j$ in
the second equality. By (\ref{3.8}) we see that $\{X_t: t\ge0\}$
solves the $({\mcr L}, {\mcr D}(\L))$-martingale problem.  \qed

The distribution of $\{X_t: t\ge0\}$ can be characterized in terms
of a dual process as follows. Let $\{M_t: t\ge 0\}$ be a
nonnegative integer-valued c\'adl\'ag Markov process with
transition intensities $\{q_{i,j}\}$ such that $q_{i,i-1} =
-q_{i,i} = i(i-1)/2$ and $q_{i,j}=0$ for all other pairs $(i,j)$.
That is, $\{M_t: t\ge 0\}$ only has downward jumps which occur at
rate $M_t(M_t-1)/2$. Such a Markov process is known as Kingman's
coalescent process. Let $\tau_0 =0$ and $\tau_{M_0} = \infty$. For
$1\le k\le M_0-1$ let $\tau_k$ denote the $k$th jump time of
$\{M_t: t\ge 0\}$. Let $\{\itGamma_k: 1\le k\le M_0-1\}$ be a
sequence of random operators which are conditionally independent
given $\{M_t: t\ge 0\}$ and satisfy
 \beqlb\label{3.11}
\bfP\{\itGamma_k = \itPhi_{ij} | M(\tau_k^-) =l\} =
\frac{1}{l(l-1)}, \qquad 1 \le i \neq j \le l,
 \eeqlb
where $\itPhi_{ij}$ is defined by (\ref{3.9}). Let $\bfC$ denote
the topological union of $\{C(\dfR^m): m=1,2, \cdots\}$ endowed
with pointwise convergence on each $C(\dfR^m)$. Let $(P_t^{(m)})
_{t\ge0}$ denote the transition semigroup of the $m$-SCBM. Then
 \beqlb\label{3.12}
Y_t = P^{(M_{\tau_k})}_{t-\tau_k} \itGamma_k
P^{(M_{\tau_{k-1}})}_{\tau_k -\tau_{k-1}} \itGamma_{k-1} \cdots
P^{(M_{\tau_1})}_{\tau_2 -\tau_1} \itGamma_1 P^{(M_0)}_{\tau_1}
Y_0, \quad \tau_k \le t < \tau_{k+1}, 0\le k\le M_0-1,
 \eeqlb
defines a Markov process $\{Y_t: t\ge0\}$ taking values from
$\bfC$. The process evolves in the time interval $[0,\tau_1)$
according to the linear semigroup $(P_t^{(M_0)})_{t\ge0}$ and then
it makes a jump given by $\itGamma_1$ at time $\tau_1$. After
that, it evolves in interval $[\tau_1,\tau_2)$ according to
$(P_t^{(M_{\tau_1})}) _{t\ge0}$ and then it makes another jump
given by $\itGamma_2$ at time $\tau_2$, and so on. Clearly,
$\{(M_t, Y_t): t\ge 0\}$ is also a Markov process. Let
$\bfE^\sigma_{m,f}$ denote the expectation related to $\{(M_t,
Y_t): t\ge0\}$ given $M_0=m$ and $Y_0=f \in C(\dfR^m)$.

\btheorem\label{t3.4} If $\{X_t: t\ge0\}$ is a continuous
$M(\dfR)$-valued process such that $\bfE \{\<1,X_t\>^m\}$ is
locally bounded in $t\ge0$ for each $m\ge1$ and $\{X_t: t\ge0\}$
solves the $({\mcr L}, {\mcr D}(\L))$-martingale problem with
$X_0=\mu$, then the distribution of $X_t$ is uniquely determined
by
 \beqlb\label{3.13}
\bfE \<f,X^m_t\> = \bfE^\sigma_{m,f} \bigg[\<Y_t, \mu^{M_t}\>
\exp\bigg\{\frac{1}{2}\int_0^t M_s(M_s-1)ds \bigg\}\bigg],
 \eeqlb
where $t\ge0$, $f \in C(\dfR^m)$ and $m\ge1$. \etheorem

{\it Proof.} It suffices to prove the equation for $Y_0=f \in
{\mcr D}^{(m)}$. In this case, we have a.s.\ $Y_t \in {\mcr
D}(\L)$ for all $t\ge0$. Set $F_\mu(m,f) = F_{m,f}(\mu) =
\<f,\mu^m\>$. By the construction (\ref{3.12}), it is not hard to
see that $\{(M_t, Y_t): t\ge 0\}$ has generator $\L^*$ given by
 \beqlb\label{3.14}
\L^* F_\mu(m,f) = F_\mu(m,G^{(m)}f) + \frac{1}{2} \sum_{i,j=1,
i\neq j}^m [F_\mu(m-1,\itPhi_{ij}f) - F_\mu(m,f)].
 \eeqlb
In view of (\ref{3.8}) and (\ref{3.14}) we have
 \beqlb\label{3.15}
\L F_{m,f}(\mu)
=
\L^* F_\mu(m,f) + \frac{1}{2} m(m-1) F_\mu(m,f).
 \eeqlb
The right hand side corresponds to a Feynman-Kac formula for the
process $\{(M_t, Y_t): t\ge 0\}$. Guided by this relation, it is
not hard to get
 \beqnn
\bfE \left[F_{m,f}(X_t)\right] = \bfE^\sigma_{m,f}
\bigg[F_\mu(M_t,Y_t) \exp\bigg\{\frac{1}{2}\int_0^t M_s(M_s-1)ds
\bigg\}\bigg],
 \eeqnn
which is just (\ref{3.13}). This formula gives in particular all
the moments of $\<f_1,X_t\>$ for $f_1\in C(\dfR)$ and hence
determines uniquely the distribution of $X_t$. We omit the details
since they are almost identical with the proofs of Dawson {\it et
al.} (2001, Theorems~2.1 and 2.2). \qed

By Theorems~\ref{t3.3} and \ref{t3.4}, the process $\{X_t:
t\ge0\}$ constructed by (\ref{3.3}) is a diffusion process. Let
$Q_t(\mu,d\nu)$ denote the distribution of $X_t$ on $M(\dfR)$
given $X_0 = \mu \in M_a(\dfR)$. The above theorem asserts that
 \beqlb\label{3.16}
\int_{M(\dfR)} \<f,\nu^m\>Q_t(\mu,d\nu) = \bfE^\sigma_{m,f}
\bigg[\<Y_t, \mu^{M_t}\> \exp\bigg\{\frac{1}{2}\int_0^t
M_s(M_s-1)ds \bigg\}\bigg]
 \eeqlb
for $t\ge0$, $m\ge1$ and $f \in C(\dfR^m)$. As in the proof of
Dawson {\it et al.} (2001, Theorem~5.1) one can extend
$Q_t(\mu,d\nu)$ to a Feller transition semigroup on $M(\dfR)$. A
Markov process on $M(\dfR)$ with transition semigroup
$(Q_t)_{t\ge0}$ given by (\ref{3.16}) is called a {\it
superprocess with coalescing spatial motion} (SCSM) with speed
$\rho$ and branching rate $\sigma(\cdot)$.

A construction of the SCSM with a general initial state $\mu \in
M(\dfR)$ is given as follows. Let $W= C([0,\infty), \dfR^+)$ and
let $\tau_0(w) = \inf\{s>0: w(s)=0\}$ for $w\in W$. Let $W_0$ be
the set of paths $w \in W$ such that $w(0)=w(t)=0$ for $t \ge
\tau_0(w)$. We endow $W$ and $W_0$ with the topology of locally
uniform convergence. Let $(q_t)_{t\ge0}$ denote the transition
semigroup of the standard Feller branching diffusion. For $t>0$
and $y>0$ let $\kappa_t(dy) = 4t^{-2} e^{-2y/t} dy$. Then
$(\kappa_t)_{t>0}$ is an entrance law for the restriction of
$(q_t)_{t\ge0}$ to $(0,\infty)$. Let $\bfQ_\kappa$ denote the
corresponding excursion law, which is the unique $\sigma$-finite
measure on $W_0$ satisfying
 \beqnn
\bfQ_\kappa\{w(t_1)\in dy_1,\cdots,w(t_n)\in dy_n\} =
\kappa_{t_1}(dy_1)q_{t_2-t_1}(y_1,dy_2)\cdots
q_{t_n-t_{n-1}}(y_{n-1},dy_n)
 \eeqnn
for $0<t_1<t_2<\cdots<t_n$ and $y_1,y_2,\cdots,y_n\in (0,\infty)$;
see, e.g., Pitman and Yor (1982) or Dawson and Li (2003, p.41) for
details. Suppose that $\{y(a,t): a\in \dfR, t\ge0\}$ is a
coalescing Brownian flow and $N(dx,dw)$ is a Poisson random
measure on $\dfR\times W_0$ with intensity $\mu(dx) \bfQ_\kappa
(dw)$. Assume that $\{y(a,t)\}$ and $\{N(dx,dw)\}$ are defined on
a standard probability space $(\itOmega,{\mcr F},\bfP)$ and are
independent of each other. As in Dawson and Li (2003), we can
enumerate the atoms of $N(dx,dw)$ into a sequence $\supp(N)
=\{(a_i,w_i): i=1,2,\cdots\}$ such that a.s.\ $\tau_0(w_{i+1}) <
\tau_0(w_i)$ for all $i\ge1$ and $\tau_0(w_i) \to 0$ as
$i\to\infty$. Let
 \beqlb\label{3.17}
\psi^\sigma_i(t) = \int_0^t \sigma(y(a_i,s)) ds
 \eeqlb
and $w_i^\sigma(t) = w_i(\psi^\sigma_i (t))$. For $t\ge0$ let
${\mcr G}_t$ be the $\sigma$-algebra generated by the family
 \beqlb\label{3.18}
\{y(a,s): 0\le s\le t, a\in \dfR\} \quad\mbox{and}\quad
\{w_i^\sigma(s): 0\le s\le t, i=1,2,\cdots\}.
 \eeqlb

\btheorem\label{t3.5} Let $X_0 = \mu$ and let
 \beqlb\label{3.19}
X_t = \sum_{i=1}^{\infty}w_i^\sigma(t)\delta_{y(a_i,t)},
\quad t>0.
 \eeqlb
Then $\{X_t: t\ge 0\}$ is a SCSM relative to $({\mcr G}_t)
_{t\ge0}$. \etheorem

\noindent{\it Proof.} For $r>0$ let $\supp_r^\sigma (N) =
\{(x_i,w_i) \in \supp(N): w^\sigma_i (r)>0\}$ and $m^\sigma(r) =
\# \{\supp_r^\sigma (N)\}$. As in Dawson and Li (2003, Lemmas~3.3
and 3.4), we have a.s.\ $m^\sigma(r) < \infty$ and there is a
permutation $\{w_{i_j}: j=1, \cdots, m^\sigma(r)\}$ of
$\supp_r^\sigma (N)$ so that $\{w_{i_j}(t): t\ge r; j=1, \cdots,
m(r)\}$ under $\bfP\{\,\cdot\, |{\mcr G}_r\}$ are independent
$\sigma$-branching diffusions which are independent of $\{x(a,t):
a\in\dfR, t\ge r\}$. By Theorem~\ref{t3.1}, $\{X_t: t\ge r\}$
under $\bfP\{\,\cdot\, |{\mcr G}_r\}$ is a Markov process with
transition semigroup $(Q_t)_{t\ge0}$. It follows that $\{X_t:
t>0\}$ is a Markov process with transition semigroup
$(Q_t)_{t\ge0}$. We shall prove that the random measure $X_t$ has
distribution $Q_t(\mu,\cdot)$ for $t>0$ so that the desired result
follows from the uniqueness of distribution of the SDSM. By
Theorem~\ref{t3.1} we can also show that
 \beqlb\label{3.20}
X_t^{(r)}
:=
\sum_{j=1}^{m(\epsilon r/\beta)} w_{k_j}(\epsilon r/\beta
+\psi_{k_j}^\sigma(t))\delta_{x_{k_j}(t)},
\quad t\ge0,
 \eeqlb
under the non-conditional probability $\bfP$ is a SDSM with
initial state
 \beqnn
X_0^{(r)} = \sum_{j=1}^{m(\epsilon r/\beta)} w_{k_j}
(\epsilon r/\beta)\delta_{a_{k_j}}.
 \eeqnn
By Shiga (1990, Theorem~3.6), $\{X_0^{(r)}: r>0\}$ is a
measure-valued branching diffusion without migration and
$X_0^{(r)} \to \mu$ a.s.\ as $r\to 0$. By the Feller property of
$(Q_t)_{t\ge0}$, the distribution of $X_t^{(r)}$ converges to
$Q_t(\mu,\cdot)$ as $r\to0$. Since $\phi_{k_j}(t)\ge \epsilon
t/\beta$, we can rewrite (\ref{3.20}) as
 \beqnn
X_t^{(r)}
:=
\sum_{j=1}^{m(\epsilon t/\beta)} w_{k_j}(\epsilon r/\beta
+\phi_{k_j}(t))\delta_{x_{k_j}(t)}.
 \eeqnn
Then for fixed $t>0$ we have $X_t^{(r)} \to X_t$ a.s.\ as $r\to0$
and hence $X_t$ has distribution $Q_t(\mu,\cdot)$. \qed

\newpage

The excursion representation (\ref{3.19}) allows us to construct
the SCSM directly without consideration of high density limits of
the corresponding coalescing-branching particle systems. This
representation also provides a useful tool for the study of the
SCSM. In particular, by (\ref{3.19}) and the proof of
Theorem~\ref{t3.4}, for each $r>0$ the process $\{X_{r+t}:
t\ge0\}$ consists of only a finite number of atoms. By this
observation and the fact $X_r \to \mu$ a.s.\ as $r\to 0$ implied
by the statements of Theorem~\ref{t3.4}, it is easy to see that
Theorems~\ref{t3.2} and \ref{t3.3} also hold for a general initial
state $\mu \in M(\dfR)$. Another application of (\ref{3.19}) is
the proof of the scaling limit theorem in the next section.

\section{A limit theorem of rescaled superprocesses}

\setcounter{equation}{0}

In this section, we show that the SCSM arises naturally as scaling
limit of the SDSM studied in Dawson {\it et al.} (2001) and Wang
(1997, 1998). In particular, the result confirms an observation of
Dawson {\it et al.} (2001) on the scaling limit of the purely
atomic SDSM.

Suppose that $h\in C(\dfR)$ is a square-integrable function with
continuous square-integrable derivative $h^\prime$. Let
$\rho(\cdot)$ be defined as in section 2. Suppose that $\sigma \in
C(\dfR)^+$ and $\inf_x \sigma(x) \ge \epsilon$ for some constant
$\epsilon>0$. We define the operator $\L$ by
 \beqlb\label{4.1}
\L F(\mu) &=& \frac{1}{2}\rho(0)\int_{\dfR}\frac{d^2}{dx^2}
\frac{\delta F(\mu)}{\delta\mu(x)}\mu(dx)  \nonumber  \\
& & +\,\frac{1}{2}\int_{\dfR^2}\rho(x-y)
\frac{d^2}{dxdy}\frac{\delta^2 F(\mu)}
{\delta\mu(x)\delta\mu(y)}\mu(dx)\mu(dy)  \nonumber \\
& & +\,\frac{1}{2} \int_{\dfR}\sigma\frac{\delta^2 F(\mu)}
{\delta\mu(x)^2}\mu(dx).
 \eeqlb
Let ${\mcr D}(\L)$ denote the collection of functions on $M(\dfR)$
of the form $F_{n,f}(\nu) := \int fd\nu^n$ with $f\in C^2(\dfR^n)$
and functions of the form
 \beqlb\label{4.2}
F_{f,\{\phi_i\}}(\nu) := f(\<\phi_1,\nu\>, \cdots, \<\phi_n,\nu\>)
 \eeqlb
with $f\in C^2(\dfR^n)$ and $\{\phi_i\}\subset C^2(\dfR)$. An
$M(\dfR)$-valued diffusion process is called a {\it superprocess
with dependent spatial motion} (SDSM) if it solves the $(\L,{\mcr
D}(\L))$-martingale problem. The existence of solution of the
$(\L,{\mcr D}(\L))$-martingale problem was proved in Dawson {\it
et al.} (2001, Theorem~5.2) and its uniqueness follows from Dawson
{\it et al.} (2001, Theorem~2.2); see also Wang (1997, 1998).

Suppose that $\sigma(x) \to \sigma_\partial$ and $\rho(x) \to 0$
as $|x| \to \infty$. Given $\theta>0$, we defined the operator
$K_\theta$ on $M(\dfR)$ by $K_\theta \mu(B) = \mu (\{\theta x:
x\in B\})$. Let $\{X_t^{(\theta)}: t\ge0\}$ be a SDSM with
parameters $(\rho,\sigma)$ and deterministic initial state
$X_0^{(\theta)} = \mu^{(\theta)} \in M(\dfR)$. Let $X^\theta_t =
\theta^{-2} K_\theta X^{(\theta)}_{\theta^2t}$ and assume
$\mu_\theta := \theta^{-2} K_\theta \mu^{(\theta)} \to \mu$ as
$\theta \to \infty$. By Dawson {\it et al.} (2001, Lemma~6.1),
$\{X^\theta_t: t\ge0\}$ is a SDSM with parameters $(\rho_\theta,
\sigma_\theta)$.

\blemma\label{l4.1} Under the above assumptions, $\{X^\theta_t:
t\ge0; \theta\ge 1\}$ is tight in $C([0,\infty), M(\bar \dfR))$.
\elemma

{\it Proof.} By Dawson and Li (2003, Theorem~3.2), $\{\<1,
X^\theta_t\>: t\ge0\}$ is a continuous positive martingale. Then
we have
 \beqnn
\bfP\bigg\{\sup_{t\ge0}\<1, X^\theta_t\> > \eta\bigg\} \le
\frac{\<1, \mu_\theta\>}{\eta}
 \eeqnn
for any $\eta>0$. That is, $\{X^\theta_t: t\ge0; \theta \ge 1\}$
satisfy the compact containment condition of Ethier and Kurtz
(1986, p.142). Let $\L_\theta$ denote the generator of
$\{X_t^\theta: t\ge0\}$ and let $F= F_{f,\{\phi_i\}}$ be given by
(\ref{4.2}) with $f\in C^2_0 (\dfR^n)$ and with each $\phi_i \in
C^2_\partial (\dfR)$ bounded away from zero. Then
 \beqnn
F(X^\theta_t) - F(X^\theta_0)
- \int_0^t \L_\theta F(X^\theta_s)ds,
\qquad t\ge0,
 \eeqnn
is a martingale and the desired tightness follows from the result
of Ethier and Kurtz (1986, p.145). \qed

Let us adopt a useful representation of the SDSM in terms of
excursions similar to the one discussed in section 3. Suppose we
have on some standard probability space $(\itOmega,{\mcr F},\bfP)$
a time-space white noise $W(ds,dy)$ on $[0,\infty) \times \dfR$
based on the Lebesgue measure and a Poisson random measure
$N_\theta (dx,dw)$ on $\dfR\times W_0$ with intensity $\mu_\theta
(dx) \bfQ_\kappa(dw)$, where $\bfQ_\kappa$ denotes the excursion
law of the standard Feller branching diffusion. Assume that
$\{W(ds,dy)\}$ and $\{N_\theta (dx,dw)\}$ are independent. We
enumerate the atoms of $N_\theta(dx,dw)$ into a sequence
$\supp(N_\theta) = \{(a_i,w_i): i=1,2, \cdots\}$ so that a.s.\
$\tau_0(w_{i+1}) < \tau_0(w_i)$ and $\tau_0(w_i) \to 0$ as $i \to
\infty$. Let $\{x^\theta(a_i,t): t\ge 0\}$ be the solution of
(\ref{2.8}) with $a_i$ replacing $a$ and $\sqrt {\theta} h_\theta
(\cdot)$ replacing $h(\cdot)$. Let
 \beqlb\label{4.3}
\psi^\theta_i(t) = \int_0^t \sigma_\theta (x^\theta(a_i,s)) ds
 \eeqlb
and $w_i^\theta(t) = w_i(\psi^\theta_i (t))$. By Dawson and Li
(2003, Theorem~3.4), the process $\{Y^\theta_t: t\ge0\}$ defined
by $Y^\theta_0 = \mu_\theta$ and
 \beqlb\label{4.4}
Y^\theta_t
=
\sum_{i=1}^{\infty} w_i^\theta(t)\delta_{x^\theta(a_i,t)},
\quad t>0.
 \eeqlb
has the same distribution on $C([0,\infty), M(\dfR))$ as
$\{X^\theta_t: t\ge0\}$. The following theorem confirms an
observation given in the introduction of Dawson {\it et al.}
(2001).

\btheorem\label{t4.1} The distribution of $\{X^\theta_t: t\ge0\}$
on $C([0,\infty), M(\dfR))$ converges as $\theta \to \infty$ to
that of a SCSM with speed $\rho(0)$, constant branching rate
$\sigma_\partial$ and initial state $\mu$. \etheorem

{\it Proof.} For any $r>0$, let $\bfQ_\kappa^r$ denote the
restriction of $\bfQ_\kappa$ to $W_r := \{w\in W_0:
\tau_0(w)>r\}$. Then we have $\bfQ_\kappa(W_r) =
\bfQ_\kappa^r(W_r) = 2/r$; see, e.g., Dawson and Li (2003). Since
$\inf_x \sigma \ge \epsilon$, we have $\psi^\theta_i(t) \ge
\epsilon t$. Then $w_i^\theta(t)=0$ for all $t\ge r$ if
$w_i(\epsilon r)=0$. Thus we only need to consider the restriction
of $N_\theta$ to $W_{\epsilon r}$ for the construction of the
process $\{Y^{\theta}_t: t\ge r\}$. To avoid triviality we assume
$\<1,\mu\>>0$. Suppose we have on a probability space the
following:

\begin{enumerate}

\item[(i)] a family of Poisson random variables $\eta_\theta$ with
parameter $\<1, \mu_\theta\> \<1, \bfQ^{\epsilon r}_\kappa\>$ such
that $\eta_\theta \to \eta$ a.s.\ as $\theta\to \infty$, where
$\eta$ is a Poisson random variable with parameter $\<1,\mu\>^{-1}
\<1, \bfQ^{\epsilon r}_\kappa\>$.

\item[(ii)] sequences of i.i.d.\ real random variables $\{a_{\theta,1},
a_{\theta,2}, \cdots\}$ with distributions $\<1, \mu_\theta\>
^{-1}$ $\mu_{\theta}(dx)$ such that $a_{\theta,i} \to a_i$ a.s.\
as $\theta \to \infty$, where $\{a_1, a_2, \cdots\}$ are i.i.d.\
real random variables with distribution $\<1, \mu\>^{-1} \mu(dx)$.

\item[(iii)] a sequence of i.i.d.\ random variables $\{\xi_1, \xi_2,
\cdots\}$ taking values from $W_{\epsilon r}$ with distribution
$\<1, \bfQ^{\epsilon r}_\kappa\> ^{-1} \bfQ^{\epsilon
r}_\kappa(dw)$.

\end{enumerate}

Under those assumptions, it is not hard to see that $\sum_{i=1}
^{\eta_\theta} \delta_{(a_{\theta,i},\xi_i)}$ and $\sum_{i=1}
^{\eta} \delta_{(a_i,\xi_i)}$ are Poisson random measures with
intensities $\mu_\theta(dx) \bfQ^{\epsilon r}_\kappa(dw)$ and
$\mu(dx) \bfQ^{\epsilon r}_\kappa(dw)$ respectively. Let
$\{x^\theta(a_{\theta,i},t): t\ge 0\}$ be the solution of
(\ref{2.8}) with $a_{\theta,i}$ replacing $a$ and $\sqrt {\theta}
h_\theta (\cdot)$ replacing $h(\cdot)$. Let
 \beqlb\label{4.5}
\psi_{\theta,i}(t) = \int_0^t \sigma_\theta
(x^\theta(a_{\theta,i},s)) ds
 \eeqlb
and $\xi_{\theta,i}(t) = \xi_i(\psi_{\theta,i}(t))$. In view of
(\ref{4.4}), the process
 \beqlb\label{4.6}
Z^{\theta}_t
:=
\sum_{i=1}^{\eta_\theta} \xi_{\theta,i}(t)
\delta_{x^\theta(a_{\theta,i},t)},
\quad t\ge r,
 \eeqlb
has the same distribution on $C([r,\infty), M(\dfR))$ as
$\{Y^{\theta}_t: t\ge r\}$ and $\{X^{\theta}_t: t\ge r\}$. By
Theorem~\ref{t2.3} it is easy to show that $\{Z^\theta_t: t\ge
r\}$ converges in distribution to
 \beqlb\label{4.7}
X_t
:=
\sum_{i=1}^{\eta} \xi_i(\sigma_\partial t)\delta_{x(a_i,t)},
\quad t\ge r,
 \eeqlb
where $\{x(a_i,t)\}$ is a system of coalescing Brownian motions.
By Theorem~\ref{t3.5}, $\{X_t: t\ge r\}$ has the same distribution
on $C([r,\infty), M(\dfR))$ as the SCSM described in the theorem.
Then the above arguments show that the distribution of
$\{X^\theta_t: t\ge r\}$ on $C([r,\infty), M(\dfR))$ converges as
$\theta \to \infty$ to that of the SCSM. The convergence is
certainly true if we consider the distributions on $C([r,\infty),
M(\bar \dfR))$. By Lemma~\ref{l4.1} it is easy to conclude that
the distribution of $\{X^\theta_t: t\ge0; \theta\ge1\}$ on
$C([0,\infty), M(\bar\dfR))$ converges to that of the SCSM. Since
all the distributions are supported by $C([0,\infty), M(\dfR))$,
the desired result follows. \qed

{\bf Acknowledgements.} We thank S.N.\,Evans and T.G.\,Kurtz for
enlightening comments on absorbing and coalescing Brownian
motions. We are indebted to a referee for a list of comments and
suggestions which helped us in improving the presentation of the
results. We are also grateful to H.\,Wang and J.\,Xiong for
helpful discussions on the subject. Dawson and Zhou were supported
by NSERC Grants, and Li was supported by NSFC Grants.

\end{document}